\documentclass[11pt,a4paper,final]{article}
\pdfoutput=1
\usepackage[utf8]{inputenc}
\usepackage[T1]{fontenc}
\usepackage[english]{babel}
\usepackage{lmodern,microtype}
\usepackage{amsmath,amssymb,amsthm,mathtools} 
\usepackage{cite} 
\usepackage{eucal,enumitem}
\usepackage[b]{esvect}

\usepackage[nottoc]{tocbibind} 
\usepackage{tikz,mdframed} 
\usepackage[colorlinks,
  linkcolor=black!40!red,%
  citecolor=green!40!black,%
  urlcolor=black,%
]{hyperref} 

\usepackage[left=2cm,right=3.5cm,top=1.5cm,bottom=1.6cm]{geometry}
\usepackage{changepage} 
\usepackage{marginnote} 
\usepackage{zref-savepos,zref-abspage} 
\usepackage{expl3} 

\usepackage{setspace}

\usepackage[notref,notcite,final%
]{showkeys} 

\usepackage[mathlines]{lineno}
\usepackage{etoolbox} 

\newcommand*\linenomathpatch[1]{%
   \expandafter\pretocmd\csname #1\endcsname {\linenomath}{}{}%
   \expandafter\pretocmd\csname #1*\endcsname{\linenomath}{}{}%
   \expandafter\apptocmd\csname end#1\endcsname {\endlinenomath}{}{}%
   \expandafter\apptocmd\csname end#1*\endcsname{\endlinenomath}{}{}%
 }
\newcommand*\linenomathpatchAMS[1]{%
    \expandafter\pretocmd\csname #1\endcsname {\linenomathAMS}{}{}%
    \expandafter\pretocmd\csname #1*\endcsname{\linenomathAMS}{}{}%
    \expandafter\apptocmd\csname end#1\endcsname {\endlinenomath}{}{}%
    \expandafter\apptocmd\csname end#1*\endcsname{\endlinenomath}{}{}%
}

\expandafter\ifx\linenomath\linenomathWithnumbers
\let\linenomathAMS\linenomathWithnumbers
\patchcmd\linenomathAMS{\advance\postdisplaypenalty\linenopenalty}{}{}{}
\else
\let\linenomathAMS\linenomathNonumbers
\fi

\linenomathpatchAMS{gather}
\linenomathpatchAMS{multline}
\linenomathpatchAMS{align}
\linenomathpatchAMS{alignat}
\linenomathpatchAMS{flalign}
\linenomathpatch{equation}

\colorlet{tn/color/theorem}{green!60!black}
\colorlet{tn/color/claim}{gray!60!white}
\colorlet{tn/color/definition}{black}
\colorlet{tn/color/question}{orange!60!white}
\colorlet{tn/color/conjecture}{orange!60!white}
\colorlet{tn/color/defi}{red!60!black}

\theoremstyle{definition}

\newtheoremstyle{slthm}
{0pt}
{0pt}
{\normalfont}
{}
{\bfseries\footnotesize}
{.\!}
{.6em}
{\thmname{#1}\thmnumber{ #2}\thmnote{\,--#3}}
\theoremstyle{slthm}

\newlength{\mythmbar}
\newlength{\myinsep}
\newlength{\myleftinsep}
\newlength{\mylen}
\newlength{\thmskip}
\mdfsetup{skipabove=\thmskip,skipbelow=\thmskip,usetwoside=false,nobreak}

\newcommand{\tndefmdstyle}[2]{%
  \mdfdefinestyle{tn/#1}{%
    linewidth=\mythmbar,%
    linecolor=tn/color/#2,%
    bottomline=false,%
    topline=false,%
    innerleftmargin=\myleftinsep,%
    leftmargin=-\myinsep,%
    innerrightmargin=\myinsep,%
    innertopmargin=1pt,%
    innerbottommargin=0pt,%
    rightmargin=\myinsep,%
    skipabove=\bigskipamount,%
    skipbelow=3pt,%
    userdefinedwidth=\mylen}}

\tndefmdstyle{base}{theorem}
\tndefmdstyle{claim}{claim}
\tndefmdstyle{definition}{definition}
\tndefmdstyle{question}{question}

\newmdtheoremenv[style=tn/base]{theorem}{Theorem}

\newcommand{\tnmdenv}[3]{%
  \newmdtheoremenv[style=#1]{#2}[theorem]{#3}}

\tnmdenv{tn/base}{lemma}{Lemma}
\tnmdenv{tn/base}{proposition}{Proposition}
\tnmdenv{tn/claim}{claim}{Claim}
\tnmdenv{tn/claim}{example}{Example}
\tnmdenv{tn/base}{corollary}{Corollary}
\tnmdenv{tn/base}{fact}{Fact}
\tnmdenv{tn/question}{question}{Question}
\tnmdenv{tn/question}{problem}{Problem}
\tnmdenv{tn/question}{conjecture}{Conjecture}
\tnmdenv{tn/definition}{definition}{Definition}
\tnmdenv{tn/definition}{remark}{Remark}
\tnmdenv{tn/definition}{hole}{Hole}

\newtheoremstyle{tnUnnumb}
{0pt}
{0pt}
{\normalfont}
{}
{\bfseries\footnotesize}
{.}
{.5em}
{\thmnote{#3}}

\theoremstyle{tnUnnumb}
\newmdtheoremenv[style=tn/base]{unnumtheorem}{Theorem}
\newmdtheoremenv[style=tn/base]{unnumconjecture}{Conjecture}

\newcommand{\defistyle}[1]{\textcolor{definition}{\textsl{#1}}}
\newcommand{\defi}[1]{%
  \sidepar{{\tiny#1}}%
  {\defistyle{#1}}}

\ExplSyntaxOn\makeatletter
\int_new:N \g_tassio_note_int
\seq_new:N\g_tassio_note_seq
\renewcommand{\defi}{\@dblarg\tn@defi} 
\renewcommand{\defistyle}[1]{%
  \textcolor{tn/color/defi}{\textsl{#1}}%
}
\def\tn@defi[#1]#2{%
  \int_gincr:N \g_tassio_note_int
  \bool_if:nTF
  {
   \int_compare_p:n
    {
     \zposy{tassionotepos\int_eval:n{\g_tassio_note_int}}
     =
     \zposy{tassionotepos\int_eval:n{\g_tassio_note_int +1}}
    }
   &&
   \int_compare_p:n
   {
    \zref@extractdefault { tassionotepage\int_eval:n{\g_tassio_note_int } }{abspage}{-1}
    =
    \zref@extractdefault { tassionotepage\int_eval:n{\g_tassio_note_int +1} }{abspage}{-1}
   }
  }
   {
    \seq_gput_right:Nn \g_tassio_note_seq {#1}
   }
   {
    \seq_gput_right:Nn \g_tassio_note_seq {#1}
    \marginnote{\fontsize{7pt}{5pt}\selectfont%
      \baselineskip=.6\baselineskip 
      \lineskip=-1.2pt 
      \seq_use:Nn \g_tassio_note_seq {,\ }}
    \seq_gclear:N \g_tassio_note_seq
   }
  \zref@label {tassionotepage\int_use:N\g_tassio_note_int}\zsaveposy {tassionotepos\int_use:N\g_tassio_note_int}
  \defistyle{#2}}

\ExplSyntaxOff

\def\tn@defi[#1]#2{{\defistyle{#2}}}

\makeatother

\makeatletter
\renewcommand{\@biblabel}[1]{\textcolor{gray}{[\,}#1\textcolor{gray}{\,]}}
  \renewcommand\@openbib@code{
    \setlength\labelwidth{2cm}%
    \setlength{\itemindent}{0cm}%
    \setlength{\leftmargin}{0cm}%
    \setlength\labelsep{.8em}%
    }
\makeatother

\newcommand{\eps}{\varepsilon}

\newcommand{\deq}{\coloneqq}

\newcommand{\sm}{\setminus}

\newcommand{\calc}{{\mathcal{C}}}
\newcommand{\cald}{{\mathcal{D}}}
\newcommand{\cale}{{\mathcal{E}}}

\newcommand{\calg}{{\mathcal{G}}}
\newcommand{\calh}{{\mathcal{H}}}

\newcommand{\cals}{{\mathcal{S}}}
\newcommand{\calt}{{\mathcal{T}}}

\newcommand{\calv}{{\mathcal{V}}}

\newcommand{\bbn}{{\mathbb{N}}}

\newcommand{\ee}{\mathrm{e}}

\newcommand{\afigure}[3]{%
  \begin{figure}
    \begin{center}
      \includegraphics{#1}
      \caption{#2}\label{#3}
    \end{center}
  \end{figure}
}

\newcommand{\digraph}[1]{\vv{#1}}

\newcommand{\ttt}{\ensuremath{\digraph{\mathrm{TT}}_3}}
\newcommand{\arr}{\to}
\newcommand{\narr}{\not\to}

\let\emptyset\varnothing
\newcommand{\Carr}{{\ensuremath{\vec C}}}
\newcommand{\Darr}{{\ensuremath{\vec D}}}
\newcommand{\Farr}{{\ensuremath{\vec F}}}
\newcommand{\Garr}{{\ensuremath{\vec G}}}
\newcommand{\Harr}{{\ensuremath{\vec H}}}

\newcommand{\Karr}{{\ensuremath{\vec K}}}

\newcommand{\Rarr}{{\ensuremath{\vec R}}}

\newcommand{\Tarr}{{\ensuremath{\vec T}}}

\newcommand*\ovl[1]{\overline{#1}}

\def\({\left(}
\def\){\right)}
\let\:=\colon

\let\geq\geqslant
\let\leq\leqslant
\let\ge\geqslant
\let\le\leqslant
\newcommand{\nrms}{\not\to}

\newcommand{\oramsey}{\ensuremath{\vec R}\,}

\newcommand{\sseq}{\subseteq}
\newcommand{\abs}[1]{|#1|}
\newcommand{\card}[1]{|#1|}
\newcommand{\set}[1]{\left\{#1\right\}}

\DeclareMathOperator{\PP}{\mathbb{P}}

\let\epsilon\varepsilon

\begin{document} 

\makeatletter
\newcommand\footnoteref[1]{\protected@xdef\@thefnmark{\ref{#1}}\@footnotemark}
\makeatother

\begin{center} 
  {\Large\noindent\textbf{Directed graphs with lower orientation
      Ramsey thresholds}\footnote{This study was financed in part by
      the Coordena\c c\~ao de Aperfei\c coamento de Pessoal de N\'ivel
      Superior, Brasil (CAPES), Finance Code 001, by FAPESP
      (2018/04876-1, 2019/13364-7) and by CNPq (423833/2018-9,
      428385/2018-4)}}

  \medskip

  Gabriel Ferreira Barros\footnote{\label{foot:IME}%
    Instituto de Matem\'atica e Estat\'{\i}stica, Universidade de
    S\~ao Paulo, Rua do Mat\~ao 1010, 05508-090 S\~ao Paulo, Brazil,
    \texttt{\{yoshi, mota\}@ime.usp.br}}\,
  \quad Bruno Pasqualotto Cavalar
    \footnote{%
      Department of Computer Science, University of Warwick, Coventry
      CV4~7AL, UK, \texttt{Bruno.Pasqualotto-Cavalar@warwick.ac.uk}}\,%
    \footnote{Supported by FAPESP (2018/05557-7)}
  \quad Yoshiharu Kohayakawa\footnoteref{foot:IME}\,\footnote{%
    Partially supported by CNPq (311412/2018-1 and 315258/2023-3)}
  \quad Guilherme Oliveira
    Mota\footnoteref{foot:IME}\,\footnote{%
    Partially supported by CNPq (306620/2020-0)}
  \quad T\'assio Naia\footnote{Centre de Reserca Matem\`{a}tica,
    Edifici C, Campus Bellaterra, 08193 Bellaterra, Spain
  \texttt{tnaia@member.fsf.org}}\,\footnote{Supported by
    FAPESP (2019/04375-5, 2020/16570-4) and by the Spanish State
    Research Agency (through the María de Maeztu Program for
    Centers and Units of Excellence
    in R\&D~CEX2020-001084-M)}
\end{center} 

\begin{abstract} 
  We investigate the threshold~$p_{\Harr}=p_{\Harr}(n)$ for the
  Ramsey-type property $G(n,p)\arr \Harr$, where $G(n,p)$ is the
  binomial random graph and $G\arr\Harr$ indicates that every
  orientation of the graph~$G$ contains the oriented graph~$\vec H$ as
  a subdigraph.  Similarly to the classical Ramsey setting,
  the upper bound $p_{\Harr}\leq Cn^{-1/m_2(\Harr)}$
  is known to hold for some constant~$C=C(\Harr)$,
  where $m_2(\Harr)$ denotes the maximum $2$-density of the underlying
  graph~$H$ of~$\Harr$.  While this upper bound is indeed the
  threshold for some~$\Harr$, this is not always the case.  We obtain
  examples arising from rooted products of orientations of sparse
  graphs (such as forests, cycles and, more generally, subcubic
  $\{K_3,K_{3,3}\}$-free graphs) and arbitrarily rooted transitive
  triangles.
\end{abstract} 

\begin{flushright}
  Dedicated to the memory of Gabriel Ferreira Barros and\\
  to Professor Jayme Luiz Szwarcfiter
  on the occasion of his 80$^{\mathrm{th}}$ birthday.
\end{flushright}

\section{Introduction}
Given a graph $G$ and an oriented graph \Harr, we write
$G\arr \Harr$ to mean that every orientation of~$G$ contains a copy
of~$\Harr$.
The \defi{orientation Ramsey number} $\oramsey(\Harr)\deq
\inf \bigl\{\, n : K_n\arr \Harr \,\bigr\}$
has been investigated by many authors~(see,
e.g.,~\cite{
  LinSaksSos83,Thomason1986,HagThom91,%
  KMO10:sumner_exact,
  MyNa2018:unavoidable,DrossHavet18}).
It was noted by Stearns (see, e.g.,\cite{ErdosMo64})
that a tournament with $2^t$~vertices
must contain a transitive tournament of order~$t$.
Therefore,
$\oramsey(\Harr)<\infty$ if and only if $\Harr$
contains no directed cycle (the ``only if'' part of the
statement follows by considering transitive tournaments).

We study the property $G\to\Harr$ in the context of the \defi{binomial
  random graph}~$G(n,p)$, which is formed from the complete
graph~$K_n$ by deleting each edge independently of all others with
probability~$1-p$.  Our goal is to find, given an acyclically oriented
graph~$\Harr$, the \defi{threshold} function $p_{\Harr}=p_{\Harr}(n)$
for the property~$G(n,p)\arr\Harr$; that is, a function satisfying
\begin{equation*}
\lim_{n\to\infty}  \PP(G(n,p)\arr \Harr) =
\begin{cases}
  0 & \text{if $p\ll p_{\Harr}$,}\\
  1 & \text{if $p\gg p_{\Harr}$,}
  \end{cases}
\end{equation*}
where $a\ll b$ (or, equivalently, $b\gg a$) means
$\lim_{n\to\infty\,} a(n)/b(n)=0$ (we speak of `the
threshold~$p_{\Harr}$', since $p_{\Harr}$ is unique up to
constant factors). We shall say that an event~$\cale$ holds
asymptotically almost surely (a.a.s.)~if $\cale$ occurs
with probability tending to~$1$ as~$n\to \infty$.

Thresholds for Ramsey-type
properties have been widely studied (see,
e.g.,~\cite{janson00:_random_graph,NPSS17:ramsey_fram} and the many
references therein).
If~$\Harr$ is acyclic, then
the property $G(n,p)\arr \Harr$ is non-trivial and monotone, and
thus has a threshold~\cite{bollobas87:_thres_funct}.

Let~$H$ be a graph.  As usual, let~$v(H)$ and $e(H)$ denote the number
of vertices and edges in~$H$, respectively.  For a graph~$H$ with~$v(H)\geq3$,
let~$\rho_2(H)=(e(H)-1)/(v(H)-2)$.  Also, let $\rho_2(K_1)=\rho_2(2K_1)=0$ and
$\rho_2(K_2)=1/2$.  An important parameter for estimating~$p_{\Harr}$ is
the \defi{maximum $2$-density}~$m_2(H)$ of $H$, which is given by
$m_2(H)=\max\{\rho_2(J):J\subset H\}$.  We also consider the
\defi{maximum density}~$m(H)$ of~$H$, given by
$m(H)=\max\left\{{e(J)}/{v(J)} : J\subset H \right\}$.  Analogous
definitions are used if $\Harr$~is an oriented graph: we denote
by~$H$ the (undirected) graph we obtain from $\Harr$ by ignoring the
orientation of its arcs, and let~$m_2(\Harr)=m_2(H)$ and
$m(\Harr) = m(H)$.

The following result gives an upper bound for~$p_{\Harr}$ for any
acyclic~$\Harr$.

\begin{theorem}\label{t:one-statement}
  For every acyclically oriented graph $\Harr$ with $m_2(\Harr)>1/2$,
  there exists a constant $C = C(\Harr)$ such that if $p \ge Cn^{-1/m_2(\Harr)}$,
  then a.a.s.\ $G(n,p) \arr \Harr$.
\end{theorem}

Theorem~\ref{t:one-statement} can be proved using the regularity
method (it suffices to combine ideas from
\cite[Section~8.5]{janson00:_random_graph} and, say,
\cite{conlon14:_klr}).  We note that, when $\Harr$ is a
transitive tournament of order at least four,
Theorem~\ref{t:one-statement} is implied by a theorem of R\"odl
and Ruci\'nski~\cite{Rodl95:_thres_ramsey}, but for other
oriented graphs this implication is no longer clear.
It will be convenient to have a variant
of Theorem~\ref{t:one-statement}, namely,
Lemma~\ref{lem:one-statement_exp} (which is a version of
Theorem~\ref{t:one-statement}  stating exponentially small
bounds on the failure probability).
Let us also remark that a
version of Theorem~\ref{t:one-statement} also appeared in a
preprint of the second author, with an almost identical proof
as the one presented here (the theorem has not been published
elsewhere).
For completeness, we note that the condition $m_2(H) > 1/2$ is
equivalent to $\Delta(H) > 1$.

Readers familiar with the so called `random Ramsey theorem' of R\"odl
and Ruci\'nski~\cite{RodlRucinski93,Rodl95:_thres_ramsey} will find it
natural that the edge probability $p=p(n)=n^{-1/m_2(\Harr)}$ appears
in Theorem~\ref{t:one-statement}, and may even guess that the
threshold~$p_{\Harr}$ is~$n^{-1/m_2(\Harr)}$ as long as~$H$ is not a
forest of stars, as this is the case in the classical Ramsey set-up involving
colourings.  This was confirmed in~\cite{Barros20:_0_statement} for
cycles and complete graphs, with a single exception.

\begin{theorem}[\cite{Barros20:_0_statement}]
    \label{t:thresholds}
    If $\Harr_t$ is an acyclic orientation of a clique or of a
    cycle with $t$ vertices, then
    \begin{align*}
      p_{\Harr_t} & =
      \begin{cases}
        \displaystyle n^{-1/m(K_4)}      &  \text{if $t=3$,}  \\
        \displaystyle n^{-1/m_2(H_t)}   & \text{if $t\ge 4$.}
      \end{cases}
    \end{align*}
\end{theorem}

Let \defi{$\ttt$} be the transitive tournament on~$3$ vertices.
Theorem~\ref{t:thresholds} with~$t=3$ tells us that
$p_{\ttt}=n^{-1/m(K_4)}=n^{-2/3}\ll n^{-1/2}=n^{-1/m_2(\ttt)}$.
Therefore, the heuristic mentioned above based on the classical
colouring case fails for~$\ttt$.  This phenomenon allows us to give a
family of oriented graphs $\Harr$
(which are not forests of stars) for
which~$p_{\Harr}\ll n^{-1/m_2(\Harr)}$.  In order to define this
family, we consider \defi{rooted oriented graphs}~$\Harr$, that is,
oriented graphs~$\Harr$ with a distinguished vertex~$r$, called the
\defi{root}.
Given an oriented graph~$\Farr$, we denote by $\Farr\circ\Harr$ the \defi{rooted
  product} of~$\Farr$ and~$\Harr$, defined as~$\Farr\circ\Harr=(V,E)$ where
\begin{align*}
  V &   =   V(\Farr)\times V(\Harr),\text{ and}\\
  E &   =   \bigl\{\bigl((f,r),(f',r)\bigr):(f,f')\in E(\Farr)\bigr\}\,
      \cup\!\!\! \bigcup_{x\in V(\Farr)} \!\!\bigl\{\bigl((x,h),(x,h')\bigr):(h,h')\in E(\Harr)\bigr\}.
\end{align*}
See Figure~\ref{fig:product} for an example of a rooted product.
Rooted products have been considered before for undirected
graphs~\cite{2009:Barriere}.

In what follows, we shall consider oriented graphs of the
form~$F\circ\Harr$ with~$\Harr$ a rooted~$\ttt$.  For brevity,
$\ttt'$~below stands for any rooted~$\ttt$.  We prove that the rooted
product $\Farr\circ\ttt'$ for some oriented graphs~$\Farr$ is such
that $p_{\Farr\circ\ttt'}\ll n^{-1/m_2(\Farr\circ\ttt')}$.

\afigure{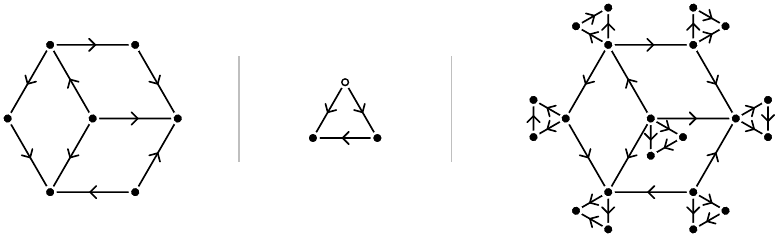}{%
  A triangle-free subcubic graph~$\Farr$,
  a rooted~\protect\ttt\ and
  their rooted product~$\Farr\circ\protect\ttt$.}{%
  fig:product}

\begin{theorem}\label{t:small-xmas}
  Let $\Farr$ be an acyclically oriented graph with $1<m_2(\Farr)<2$.
  Then
  \begin{equation}
    \label{eq:2}
    p_{\Farr\circ\ttt'}\ll n^{-1/m_2(\Farr\circ\ttt')}.
  \end{equation}
\end{theorem}

One can check that if~$F$ is a $\{K_3,K_{3,3}\}$-free graph with
maximum degree at most three that is not a forest, then any acyclic
orientation~$\Farr$ of~$F$ satisfies the hypothesis of
Theorem~\ref{t:small-xmas} (see
Proposition~\ref{p:m2<2}\,\ref{i:F_b}), and hence~\eqref{eq:2} holds.
While the hypothesis of Theorem~\ref{t:small-xmas} requires that~$F$
should contain a cycle ($m_2(\Farr)>1$ is required), the conclusion of
Theorem~\ref{t:small-xmas} also holds when~$F$ is a forest: to see
this, it suffices to take~$\delta<1/6$ in the result below.  In fact,
the result below is more general in the sense that it allows us to
consider forests~$F$ with growing order.

\begin{theorem}\label{t:large-xmas-tree}
  Let $p=n^{\delta-2/3}$ where $2/21 < \delta \le 1/6$ is a constant.
  Then a.a.s.\ $G(n, p)\arr\Farr\circ\ttt'$ for every oriented
  forest~$\Farr$ with $v(\Farr)\le bn^{7\delta-2/3}/(\log n)^2$,
  where~$b$ is some positive absolute constant.
\end{theorem}

Theorem~\ref{t:large-xmas-tree} guarantees the presence of
somewhat large oriented trees in every orientation of~$G(n,p)$.
If~$\delta<1/6$, then
for any oriented forest~$\Farr$ we have
$n^{\delta-2/3}\ll n^{-1/2}=n^{-1/m_2(\Farr\circ\ttt')}$, and indeed
Theorem~\ref{t:large-xmas-tree} extends Theorem~\ref{t:small-xmas}.
Also, if~$\delta>2/21$, then $7\delta-2/3>0$ and the forest~$\Farr$ of
Theorem~\ref{t:large-xmas-tree} may be chosen to have more
than~$n^c$ vertices for some constant~$c>0$.
And if~$\delta=1/6$, then
$v(\Farr\circ\ttt')=\Omega(n^{1/2}/(\log n)^2)$.

Lemma~\ref{lem:one-statement_exp}
(a variant of Theorem~\ref{t:one-statement})
and Theorem~\ref{t:small-xmas}
are proved in Section~\ref{s:small-xmas}
using lemmas from~Section~\ref{s:containers}.
Theorem~\ref{t:large-xmas-tree} is proved in
Section~\ref{s:large-xmas-tree}.

\smallskip

Given a set~$V$ and a positive integer~$\ell$, we denote by
$\binom{V}{\ell}$ the collection of $\ell$-element subsets of $V$\!,
by $2^V$ the collection of subsets of~$V$\!, and by $[\ell]$ the
set~$\{1,2,\dots,\ell\}$.  To avoid uninteresting technicalities, we
omit floor and ceiling signs whenever they are not important. All
unqualified logarithms are base $\ee$ (where $\ee$ denotes Euler's
constant).

\smallskip

The main results of this work were announced in the extended
abstract~\cite{barros21:_lower_orien_thres}.

\section{A container lemma}
\label{s:containers}

Let~$\Harr$ be an acyclically oriented graph.  In this section, we derive
a container lemma for graphs~$G$ that admit $\Harr$-free orientations
from the celebrated container method of Balogh, Morris and
Samotij~\cite{Balogh15} and Saxton and Thomason~\cite{Saxton15}.

\subsection{Preliminaries}
\label{sec:preliminaries}

We begin with a simple saturation result, somewhat in the spirit
of~\cite{1983:Erdos_Simonovits_supersaturation}.

\begin{lemma}
    \label{l:dir_rms_quant}
    Let $\Harr$ be an acyclically oriented graph on $h$ vertices,
    and let $R := \Rarr(\Harr )$.
    The following holds for every $n \geq R$.
    For every $F \sseq E(K_n)$,
    if there exists an orientation $\Farr$ of $F$ such that~$\Farr$ has at most
    $\bigl(2\binom{R}{h}\bigr)^{-1} \binom{n}{h}$
    copies of $\Harr$, then
    \begin{equation*}
        \bigl| E(K_n) \sm F\bigr|
        \geq
        (2R^2)^{-1}\,n^2.
    \end{equation*}
\end{lemma}
\begin{proof}
    For convenience,
    let
    $\eps := \bigl(2\binom{R}{h}\bigr)^{-1}$.
    Let $F \sseq E(K_n)$
    be such that there exists an orientation
    $\Farr$
    of
    $F$
    with at most
    $\eps \binom{n}{h}$
    copies of
    $\Harr$.
    Let $\Karr$ be an orientation of $K_n$ which agrees
    with the orientation $\Farr$ of $F$.
    Let
    \begin{equation*}
        \cals :=
        \biggl\{S \in \binom{V(\Karr)}{R} :
    E(\Karr[S]) \sseq \Farr\biggr\}.
    \end{equation*}
        That is,
        the family
        $\cals$
        is the collection of
        all
        $R$-element subsets $S$ of $V(\Karr)$
        such that
        every arc of $\Karr[S]$
        is contained in $\Farr$.
    By definition of $R$,
    every
    $R$-element
    subset of
    the vertices of
    $\Karr$
    contains
    at least one
    copy
    of $\Harr$.
    This means that, for every $S \in \cals$,
    there exists one copy of $\Harr$
    in
    $E(\Karr[S])$.
    Moreover,
    every copy of
    $\Harr$ in $\Karr$
    is contained in at most
    $\binom{n-h}{R-h}$
    $R$-element subsets.
    Double-counting the pairs
    $(S, \Harr')$
    where
    $S \in \cals$
    and
    $\Harr'$ is a copy of $\Harr$
    contained in $S$
    yields
    \begin{equation*}
        \abs{\cals}
        \leq
        \eps
        \binom{n-h}{R-h}
        \binom{n}{h}
        =
        \frac{1}{2}
        \frac{\binom{n-h}{R-h}}{\binom{R}{h}}
        \binom{n}{h}
        =
        \frac{1}{2}
        \binom{n}{R}.
    \end{equation*}
    This implies that the set
    $\ovl{\cals}$ defined as $\ovl{\cals}
        :=
        \binom{V(\Karr)}{R} \sm
        \cals$
    satisfies
    $\abs{\ovl{\cals}} \geq (1/2)\binom{n}{R}$.
    Every set
    $S \in
    \ovl{\cals}$
    induces at least one arc
    $e \in E(\Karr) \sm \Farr$.
    Moreover, every arc
    $e \in E(\Karr) \sm \Farr$
    is contained in at most
    $\binom{n-2}{R-2}$
    $R$-element subsets.
    Double-counting the pairs
    $(S,e)$
    where
    $S \in
    \ovl{\cals}$
    and
    $e \in E(\Karr[S])\setminus\Farr$
    we obtain
    \begin{equation*}
        \bigl| E(\Karr) \sm \Farr \bigr|
        \geq
        \frac{\abs{\ovl{\cals}}}{\binom{n-2}{R-2}}
        \geq
        \frac{1}{2}
        \frac {\binom{n}{R}}
              {\binom{n-2}{R-2}}
        \geq
        \frac{1}{2R^2}
        \,
        n^2.
    \end{equation*}
    The desired result now follows by observing that
    $\bigl| E(K_n) \sm F \bigr|
    =
    \bigl| E(\Karr) \sm \Farr \bigr|$.
\end{proof}

We now recall the hypergraph container lemma.  We need to introduce
some notation first.  Let $\calh$ be an
$l$-uniform hypergraph, and let~$v \in V(\calh)$.  For each
$J \sseq V(\calh)$, we call $d(J) := \bigl|\{e \in E(\calh) : J \sseq e\}\bigr|$
the \defi{degree of $J$}, and write $d(v)$ for~$d\bigl(\set{v}\bigr)$.  For each
$j \in [l]$, the \defi{maximum $j$-degree} of~$v$ is
$d^{(j)}(v) := \max\set{d(J) : v \in J \in \binom{V(\calh)}{j}}$.  We
also let
\begin{equation*}
    d_j :=
    \frac{1}{v(\calh)}
    \sum_{v \in V(\calh)}d^{(j)}(v),
\end{equation*}
and note
that
$d_1$
is the average degree of $\calh$.
Finally, for $\tau > 0$, let
\begin{equation*}
    \delta_j := \frac{d_j}{d_1 \tau^{j-1}}
\end{equation*}
and define the \defi{co-degree function}
$\delta(\calh,\tau)$
by
\begin{equation}
  \label{eq:delta_def}
    \delta(\calh,\tau) :=
    2^{\binom{l}{2}-1}
    \sum_{j=2}^l
    {2^{-\binom{j-1}{2}}}
    \delta_j,
\end{equation}
where $\binom{1}{2}=0$.  We shall apply the following hypergraph
container lemma, given in~\cite{2021:Han_et_al_ramsey-type-numbers}.

\begin{theorem}[{\cite[Theorem 2.1]{2021:Han_et_al_ramsey-type-numbers}}]
    \label{thm:container}
    Let $0< \eps,\, \tau < 1/2$ and~$l\geq2$ be given.
    There exist integers
    $K=K(l)$ and $s=s(l,\eps)$ such that
    the following holds.
    Let
    $\calh = (V,E)$
    be an $l$-uniform hypergraph and suppose~$\tau$
    is such that
    $\delta(\calh,\tau) \leq \eps/(12l!)$.
    Then, for every independent set
    $I \sseq V$ in $\calh$,
    there exist an $s$-tuple
    $T = (T_1,\dots,T_s)$
    of subsets of~$V$
    and a subset
    $C = C(T) \sseq V$ depending only on $T$
    such that
    \begin{enumerate}[label=(\alph*)]
        \item\label{i:container/a}
            $\bigcup_{i \in [s]} T_i \sseq I \sseq C$,
        \item\label{i:container/b}
            $e(\calh[C]) \leq \eps e(\calh)$, and
        \item\label{i:container/c}
            $\abs{T_i} \leq K\tau\abs{V}$
            for every $i \in [s]$.
    \end{enumerate}
\end{theorem}

Theorem~\ref{thm:container} above is a version of
\cite[Corollary~3.6]{Saxton15}.  For completeness, we mention that,
in~\cite{2021:Han_et_al_ramsey-type-numbers}, explicit values are
given for the constants~$K=K(l)$ and $s=s(l,\epsilon)$: $K=800l(l!)^3$
and~$s=\lfloor K\log(1/\eps)\rfloor$.  We now define the hypergraph on
which we shall apply Theorem~\ref{thm:container}.

\begin{definition}
  Let $n \in \bbn$ be given.  The \defi{complete digraph} $\Darr_n$ is
  the digraph with vertex set $[n]$ and arc
  set~$E(\Darr_n) := \bigl([n] \times [n]\bigr) \sm \bigl\{(v,v) : v
  \in [n]\bigr\}$.
\end{definition}

\begin{definition}
  \label{def:DnH}
  Let $\Harr$ be an oriented graph with $l$ arcs and let $n \in \bbn$
  be given.  The hypergraph $\cald(n, \Harr) = (\calv, \cale)$ is the
  $l$-uniform hypergraph with vertex set $\calv := E(\Darr_n)$ and
  edge set
    \begin{equation*}
        \cale :=
        \set{B \in \binom{\calv}{l} :
        \text{the arcs of $B$ form a digraph isomorphic to
        $\Harr$}
    }.
\end{equation*}
\end{definition}

We close this section estimating the quantity
$\delta(\cald(n, \Harr), \tau)$ defined in~\eqref{eq:delta_def} for a
certain relevant value of~$\tau$.

\begin{lemma}
    \label{lemma:deg_cont}
    Let $\Harr$ be an oriented graph
    with order~$h$
    and $l\ge 2$~arcs.
    Let $\tau := Dn^{-1/m_2(\Harr)}$, where $D\geq 1$ is a constant.
    We have
    \begin{equation}
      \label{eq:deg_cont}
        \delta(\cald(n, \Harr), \tau)
        \leq
        2^{\binom{l}{2}}
        h^{h-2}
        D^{-1}.
    \end{equation}
\end{lemma}

\begin{proof}
    For convenience,
    set
    $\calh :=
    \cald(n, \Harr)$.
    Note that $\calh$ is regular,
    i.e., $d(e)=d_1$ for all~$e\in E(\Harr)$.
    Given
    $J \sseq V(\calh)$ with~$J\neq\emptyset$ and $2\le |J|\le l$,
    let
    \begin{equation*}
        V_J
        :=
        \bigcup_{(a,b) \in J}
        \set{a,b} \sseq [n].
    \end{equation*}
    Note that
    $(V_J,J)$
    is the subdigraph of
    $\Darr_n$
    induced by the set
    of arcs~$J$.
    Given $J\subseteq V(\calh)$ and an edge $e_0\in J$,
    choose a hyperedge $F\in E(\calh)$ containing~$e_0$
    uniformly at random.
    By~symmetry, $V_F\setminus e_0$ is a uniformly-chosen
    $(h-2)$-subset of $[n]\setminus e_0$.  Since $J\subseteq F$
    implies $V_J \subseteq V_F$, we have
    \begin{equation}
      \label{eq:4}
      \frac{d(J)}{d(e_0)}
      =    \PP(J\subseteq F)
      \le  \PP(V_J\subseteq V_F)
      =    \binom{n-|V_J|}{h-|V_J|}\binom{n-2}{h-2}^{-1}\!\!
      \le  \left(\frac{h}{n}\right)^{|V_J|-2}.
    \end{equation}
    For all $\ell \in [l]$, let
    \[    f(\ell)
      :=
      \min\{v(\Harr'):\Harr' \sseq \Harr\text{ and
      }e(\Harr')=\ell\}.
    \]
    It follows from~\eqref{eq:4} that, for every $2\leq j\leq l$,
    we have $d^{(j)}(v)/d_1 \leq (h/n)^{f(j)-2}$ for any~$v\in
    V(\calh)$.
    Since $f(j) \leq h$, this gives us that
    \begin{equation*}
        \frac{d_j}{d_1}
        =
        \frac{1}{v(\calh)}
        \sum_{v \in V(\calh)}
        \frac{d^{(j)}(v)}{d_1}
        \leq
        \frac{1}{v(\calh)}
        \sum_{v \in V(\calh)}
        h^{f(j)-2}n^{2-f(j)}
        =
        h^{f(j)-2}n^{2-f(j)}
        \leq
        h^{h-2}n^{2-f(j)}.
    \end{equation*}
    We furthermore have that
    \begin{align}
        \label{eq:dgc1}
        \delta_j
        =
        \frac{d_j}{d_1 \tau^{j-1}}
        \leq
        h^{h-2}n^{2-f(j)}
        \tau^{1-j}
        \leq
        h^{h-2}
        n^{2-f(j)+(j-1)/m_2(\Harr)}
        D^{1-j}.
    \end{align}
    Observe now that, by definition,
    we have
    $m_2(\Harr) \geq (j-1)/(f(j)-2)$ for all $j\ge 2$.
    From this we may derive
    ${2-f(j)+(j-1)/m_2(\Harr)} \leq 0$.
    Therefore, we can conclude from~(\ref{eq:dgc1})
    that
    \begin{equation*}
        \label{eq:dgc2}
        \delta_j
        \leq
        h^{h-2}
        D^{1-j}
        \leq
        h^{h-2}
        D^{-1}.
    \end{equation*}
    We can finally bound the co-degree function
    $\delta(\calh, \tau)$ by observing that
    \begin{equation*}
        \delta(\calh, \tau)
        =
        2^{\binom{l}{2}-1}
        \sum_{j=2}^l
        {2^{-\binom{j-1}{2}}}
        \delta_j
        \leq
        2^{\binom{l}{2}-1}
        h^{h-2}D^{-1}
        \sum_{j=2}^l
        2^{-\binom{j-1}{2}}
        \leq
        2^{\binom{l}{2}}
        h^{h-2}D^{-1},
      \end{equation*}
      which establishes~\eqref{eq:deg_cont}.
\end{proof}

\subsection{A container lemma for graphs with \Harr-free orientations}

Let~$\Harr$ be an acyclically oriented graph.  Applying
Theorem~\ref{thm:container} to the hypergraph~$\cald(n,\Harr)$ from
Definition~\ref{def:DnH} gives us a container lemma for $\Harr$-free
digraphs.  We need something a little different: we need a container
lemma for graphs~$G$ that admit $\Harr$-free orientations.  This is
given in the lemma below.

\begin{lemma}\label{lem:container_dir_rms}
    Let
    $\Harr$
    be an acyclically oriented graph
    with at least two~arcs.
    There exist
    positive real numbers~$\alpha< 1$ and~$c$
    and a positive integer~$s$
    such that
    the following holds
    for every large enough~$n$.
    Let~$r=\lfloor cn^{2-1/m_2(\Harr)}\rfloor$.
    If~$G$ is a graph of order~$n$
    and
    $G \nrms \Harr$,
    then there exist
    an $s$-tuple
    $T=(T_1,\dots,T_s) \in \bigl(2^{E(G)}\bigr)^s$
    and
    a set $\calc=\calc(T)\subseteq 2^{E(K_n)}$ of
    size at most $2^r$,
    depending only on~$T$, such that
    \begin{enumerate}[label=(\alph*)]
    \item\label{i:cdr/a}
      $\bigcup_{i \in [s]} T_i \sseq E(G) \sseq C$
      for some~$C\in \calc$,
    \item
        \label{i:cdr/d}
      $\bigcup_{i \in [s]} T_i \sseq C$ for every $C\in \calc$.
    \item\label{i:cdr/b}
      $\bigl|E(K_n)\sm C\bigr| \geq \alpha n^2$ for
      every~$C\in\calc$, and
    \item\label{i:cdr/c}
      $\bigl|\bigcup_{i \in [s]}T_i\bigr| \leq r$.
    \end{enumerate}
\end{lemma}
\begin{proof}
    We apply Theorem~\ref{thm:container}.
    Suppose $\Harr$ has $h$ vertices and $l\geq2$ arcs and let
    $\alpha = (2R^2)^{-1}$, where  $R:=\Rarr(\Harr)$.
    In what follows, we assume that~$n$ is large enough for
    our inequalities to hold.
    Set
    $\calh := \cald(n, \Harr)$
    and
    $\eps := \bigl(2 \binom{R}{h}  e(\cald(h, \Harr))\bigr)^{-1}$
    and
    let
    $\tau := Dn^{-1/m_2(\Harr)}$,
    where
    \begin{equation*}
        D:=\frac{12\,l!\,2^{\binom{l}{2}}h^{h-2}}{\eps}.
    \end{equation*}
    By
    Lemma~\ref{lemma:deg_cont},
    we have
    $\delta(\calh, \tau) \leq \eps/(12l!)$.
    Theorem~\ref{thm:container}
    now
    gives us
    numbers $K=K(l)$ and $s=s(l,\epsilon)$ satisfying the conclusion
    of that theorem.

    Let $G$ be a graph
    on
    $n$ vertices
    such that
    $G \nrms \Harr$.
    There exists an orientation
    $\Garr$ of $G$ that is $\Harr$-free.
    Therefore,
    the set
    $E(\Garr)$
    is an independent set
    of $\calh$.
    Let
    $\Tarr = (\Tarr_1,\dots,\Tarr_s)$
    be an $s$-tuple
    of sets of arcs
    and
    $\Carr = \Carr(\Tarr)$ a set of arcs
    such as
    given by~Theorem~\ref{thm:container}
    applied to~$E(\Garr)$.
    For $i \in [s]$,
    let $T_i$ be the
    set of underlying (undirected) edges
    of the arcs in~$\Tarr_i$.
    Let
    \begin{equation}
      \label{eq:T_def}
      T=(T_1,\dots,T_s) \in \bigl(2^{E(K_n)} \bigr)^s.
    \end{equation}
    Define $C$ analogously from~$\Carr$, i.e., let~$C$ be the set of
    edges that we obtain by replacing each arc of~$\Carr$ by its
    underlying edge (we remark in passing that it may happen that
    $|C|<|\Carr|$, as~$\Carr$ may contain digons).  By
    Theorem~\ref{thm:container}\,\ref{i:container/a}, we have
    \begin{equation}
      \label{eq:concl_a}
      \bigcup_{i \in [s]}T_i \sseq E(G) \sseq C.
    \end{equation}
    Clearly, $e(\cald(h, \Harr))$
    is the number of copies
    of
    $\Harr$
    in any subset of~$h$ vertices
    of~$\Darr_n$,
    whence
    \begin{equation*}
        e(\calh) = \binom{n}{h}e(\cald(h, \Harr)).
    \end{equation*}
    Therefore,
    by Theorem~\ref{thm:container}\,\ref{i:container/b}
    we conclude that $\Carr$
    spans at most
    $\eps e(\calh) = \bigl(2 \binom{R}{h}\bigr)^{-1}\binom{n}{h}$
    copies of~$\Harr$.
    Let~$\Farr$ be an orientation of~$C$ such that
    $\Farr\subseteq \Carr$ (equivalently, let~$\Farr$ be obtained
    from~$\Carr$ by dropping an arbitrary arc from each digon
    in~$\Carr$).
    Note that the number of copies of~$\Harr$ in~$\Farr$
    is at most the number of copies of~$\Harr$ in~$\Carr$.
    Since the underlying graph of~$\Farr$ is~$C$,
    Lemma~\ref{l:dir_rms_quant} gives, by the choice
    of
    $\alpha$, that
    \begin{equation}
      \label{eq:concl_b}
      \abs{E(K_n)\sm C} \geq \alpha n^2.
    \end{equation}
    Finally,
    by letting
    $c := sKD$,
    we get by Theorem~\ref{thm:container}\,\ref{i:container/c} that
  \begin{equation}
      \label{eq:concl_c}
      \sum_{i \in [s]}
      \card{T_i}
      \leq
      s\max_{i \in [s]}\card{T_i}
    \leq
    \lfloor sK\tau v(\calh)\rfloor
    \le
    \bigl\lfloor cn^{2-1/m_2(\Harr)}\bigr\rfloor
    =r.
  \end{equation}

  We are supposed to produce a set
  $\calc=\calc(T)\subseteq 2^{E(K_n)}$ with $|\calc| \le 2^r$
  satisfying~\ref{i:cdr/a}, \ref{i:cdr/d}, \ref{i:cdr/b}
  and~\ref{i:cdr/c} of our lemma.
  Note that, because
  of~\eqref{eq:concl_a},
  if~$C\in\calc$,
  we are in good shape with respect to~\ref{i:cdr/a} in our
  lemma.  Owing to~\eqref{eq:concl_c}, we are also fine with respect
  to~\ref{i:cdr/c} in our lemma.
  We shall now see how to produce~$\calc$.

  We know that~$\Carr = \Carr(\Tarr)$ depends solely on~$\Tarr$, and
  we have produced~$C\subseteq E(K_n)$ from~$\Carr$.  Hence~$C$
  depends only on~$\Tarr$.  However, in the procedure above, $T$~as
  defined in~\eqref{eq:T_def} may arise from other $s$-tuples
  $\Tarr'=(\Tarr_1',\dots,\Tarr_s')$
  with~$\Carr(\Tarr')\neq\Carr(\Tarr)$, and hence we cannot say
  that~$C$ is solely determined by~$T$.  The solution is simple: we
  consider all possible $\Tarr'=(\Tarr_1',\dots,\Tarr_s')$ that give
  rise to~$T$ and for which $\Carr(\Tarr')$ is defined,
  to produce the set~$\calc$ of all the corresponding
  containers~$C'\subseteq E(K_n)$ coming from~$\Carr(\Tarr')$.  Owing
  to~\eqref{eq:concl_c}, the number of such~$\Tarr'$ is at most~$2^r$
  (each edge $\{x,y\}$ in each~$T_i$ in~$T$ can be oriented
  as~$\vv{xy}$ or~$\vv{yx}$ to produce~$\Tarr_i'$ in~$\Tarr'$
  and~$\Tarr_i'$ cannot contain both~$\vv{xy}$ and~$\vv{yx}$).  Hence
  the set~$\calc$ of the corresponding containers~$C'$ coming from all
  the~$\Carr(\Tarr')$ is such that $|\calc|\leq2^r$.  Because
  of~\eqref{eq:concl_b}, condition~\ref{i:cdr/b} of our lemma is
  satisfied.  Finally,
  any other $s$-tuple $\Tarr'$ that gives rise to $T$
  is such that
  $\bigcup_{i\in[s]}T_i = \bigcup_{i\in[s]}T_i' \sseq C'$,
  where~$C'\in\calc$ is the container produced by $\Tarr'$.  This
  implies~\ref{i:cdr/d}.
  Therefore, there is an $s$-tuple~$T$ and a set $\calc=\calc(T)$
  that depends only on~$T$ with the desired requirements.
\end{proof}

\section{Small graphs with low thresholds}
\label{s:small-xmas}

We prove Theorem~\ref{t:small-xmas} in this section.
We begin with a simple proposition.
\begin{proposition}
  The following statements hold.
  \label{p:m2<2}
  \begin{enumerate}[label=(\alph*)]
  \item\label{i:F_a} Let~$\Farr$ be an oriented graph with
    $m_2(\Farr)<2$ and let~$\ttt'$ be a rooted~$\ttt$.  Then
    \[m_2(\Farr\circ\ttt')=2.\]
  \item\label{i:F_b} Let~$F$ be a~$\{K_3,K_{3,3}\}$-free graph~$F$ with $\Delta(F)\leq3$ that is not a
    forest.  Then
    \[1<m_2(F)<2.\]
  \end{enumerate}
\end{proposition}

\begin{proof}
  Let~$J$ be a graph with $v(J)\geq3$ and
  $\rho_2(J)=(e(J)-1)/(v(J)-2)<2$ and let~$T$ be a~$K_3$ with
  $|V(J)\cap V(T)|\leq1$.  Let $J'=J\cup T'$ where $T'$ is a subgraph
  of~$T$.  One can check that $\rho_2(J')<2$ (roughly speaking,
  $J'$~can be obtained from~$J$ by adding~$h$ new vertices and at
  most~$3h/2$ edges for some~$h$, and hence $\rho_2(J')<2$ follows
  from a simple calculation).  If~$v(J)=2$, it is again true that
  $\rho_2(J')<2$.  These observations imply that
  $m_2(\Farr\circ\ttt')\leq2$ for any~$\Farr$ as in~\ref{i:F_a} of our
  proposition.  Since $m_2(\Farr\circ \ttt') \ge m_2(\ttt') = 2$,
  assertion~\ref{i:F_a} follows.

  Now let~$F$ be as in~\ref{i:F_b}.  Since~$F$ contains a cycle, we
  have $m_2(F)>1$.  Using the fact that~$F$ is $\{K_3,K_{3,3}\}$-free,
  one can check that $m_2(F)<2$ if~$v(F)\leq6$.  For any order~$t$
  graph~$J$ with $\Delta(J)\leq3$, we have $\rho_2(J)\leq f(t)$, where
  $f(t)=(3t/2-1)/(t-2)$.  It now suffices to observe that $f(6)=2$ and
  that~$f(t)$ is a strictly decreasing function.
\end{proof}

  We shall need
the following variant of Theorem~\ref{t:one-statement}.

\begin{lemma}\label{lem:one-statement_exp}
  For every acyclically oriented graph~$\Harr$ with~$m_2(\Harr)>1/2$,
  there are constants~$B$ and~$\beta>0$ such that if
  $p=p(n)\geq Bn^{-1/m_2(\Harr)}$,
  then
  \begin{equation}
    \label{eq:1-fail-prob}
    \PP(G(n,p)\not\arr\Harr)\leq \ee^{-\beta pn^2}.
  \end{equation}
\end{lemma}
\begin{proof}
  We first show that
  there exist constants $A$ and $\gamma > 0$ such that
  if~$q=q(n)=An^{-1/m_2(\Harr)}$ (note the equality) then
  \begin{equation}
    \PP(G(n,q)\not\arr\Harr)\leq \ee^{-\gamma qn^2}.
    \label{e:exposition}
  \end{equation}

  Let~$\Harr$ as in the statement of our lemma be given.  Let
  $\alpha$, $c$ and~$s$ be as given by
  Lemma~\ref{lem:container_dir_rms} for~$\Harr$.  Examining the
  statement of
  Lemma~\ref{lem:container_dir_rms}, one sees that
  we may suppose that~$c\geq \ee2^{s+2}$.  Let~$A$ be such that
  \begin{equation}
    \label{eq:A_cond}
    \frac{\log A}{A}\leq\frac{\alpha}{2c}
    \quad
    \text{and}
    \quad
    A \geq \frac{c}{2^{s-1}},
  \end{equation}
  and let $q=q(n) = An^{-1/m_2(\Harr)}$.  We show that if
  $\gamma=\alpha/2$, then~\eqref{e:exposition} holds for every large
  enough~$n$.
  Let $r:=\bigl\lfloor cn^{2-1/m_2(\Harr)}\bigr\rfloor$ and note
  that~$r\to\infty$ as~$n\to\infty$ because~$m_2(\Harr)>1/2$.

  Let~$G$
  be a graph on~$[n]$ such that $G \nrms \Harr$.  By
  Lemma~\ref{lem:container_dir_rms}, there exist an $s$-tuple
  $T=(T_1,\dots,T_s)\in(2^{E(K_n)})^s$ and a set
  $\calc(T)\subseteq 2^{E(K_n)}$ with $|\calc(T)|\le2^r$, depending
  only on~$T$, such that Lemma~\ref{lem:container_dir_rms}\,\ref{i:cdr/a},
  \ref{i:cdr/d},
  \ref{i:cdr/b}
  and~\ref{i:cdr/c} hold.  Let~$C\in\calc(T)$ be as in
  Lemma~\ref{lem:container_dir_rms}\,\ref{i:cdr/a}.  Then
  \begin{align}
    \label{eq:prob1}
    \bigcup_{i \in [s]} T_i
    & \sseq E(G) \sseq C
  \end{align}
  and
  \begin{align*}
    \abs{E(K_n)\sm C}
    & \geq
      \alpha n^2.
  \end{align*}
  Let $D_C := E(K_n) \sm C$.  Since $E(G) \sseq C$, we have
  \begin{align}
    \label{eq:prob2}
    E(G) \cap D_C
    & = \emptyset.
  \end{align}
  For convenience, let $\calg$ be the family of all graphs $G$
  on~$[n]$ such that $G \nrms \Harr$.  Let us summarise what we did
  above: given~$G\in\calg$, we found certain objects
  $T=(T_1,\dots,T_s)\in(2^{E(K_n)})^s$,
  $\calc(T)\subseteq 2^{E(K_n)}$, and~$C\in\calc(T)$.  Let~$\calt_n$
  be the set of all~$T$ that arise in this fashion when we consider
  all~$G$ in~$\calg$.  By construction, for each $T\in\calt_n$, we
  have a certain associated family of sets~$\calc(T)$.

  We now proceed as follows.  For $T=(T_1,\dots,T_s)\in\calt_n$, let
    \begin{align*}
      \calg_T'
      & :=
      \set{G\in\calg : T_i \sseq E(G)\text{ for all }i \in [s]}.\\
      \shortintertext{For any $C\subseteq E(K_n)$, let}
      \calg_C''
      & :=
        \set{G\in\calg : E(G) \cap D_C= \emptyset}.
    \end{align*}
    Our constructions above of~$T\in\calt_n$ and~$C\in\calc(T)$
    given~$G\in\calg$ (see, in particular, \eqref{eq:prob1}
    and~\eqref{eq:prob2}) show that
    \begin{align*}
      \calg
      \sseq
      \bigcup_{T\in\calt_n}\,\bigcup_{C\in\calc(T)}\calg_T' \cap\calg_C''.
    \end{align*}
    For any $T\in\calt_n$, the sets $T_i$ occurring in~$T$ and
    the~$D_C$ for~$C\in\calc(T)$ are disjoint (see
    Lemma~\ref{lem:container_dir_rms}~\ref{i:cdr/d}).  Therefore, the
    events $\{G(n,q) \in \calg_T'\}$ and $\{G(n,q) \in \calg_C''\}$
    are independent for any $T\in\calt_n$ and any $C\in\calc(T)$.  We
    conclude that
    \begin{equation*}
        \PP\bigl(G(n,q) \in \calg\bigr)
        \leq
        \sum_{T \in \calt_n}\,
        \sum_{C\in\calc(T)}
        \PP\bigl( G(n,q) \in \calg_T' \bigr)
        \,\PP\bigl( G(n,q) \in \calg_C'' \bigr).
    \end{equation*}
    Since $\abs{D_C} \geq \alpha n^2$ for every $C\in\bigcup_{T\in\calt_n}\calc(T)$, we have
    \begin{align}
        \PP\bigl( G(n,q) \in \calg_C'' \bigr)
        & \leq
        (1-q)^{\alpha n^2}
        \leq
        \exp(-\alpha n^2 q).\nonumber\\
      \shortintertext{Moreover, we also have}
      \sum_{T \in \calt_n}
      \PP\bigl( G(n,q) \in \calg_T' \bigr)
        & \leq
          \sum_{T \in \calt_n}
          q^{\big|\bigcup_{i \in [s]} T_i\big|}.\nonumber\\
      \shortintertext{Since $|\calc(T)|\le 2^r$, it follows that}
      \label{eq:sb}
      \PP\bigl(G(n,q) \in \calg\bigr)
        & \leq2^r
          \exp(-\alpha n^2 q)
          \sum_{T \in \calt_n}
          q^{\big|\bigcup_{i \in [s]} T_i\big|}.
    \end{align}
    We now proceed to bound the sum in
    (\ref{eq:sb}).
    For every integer $k$
    with
    $0 \leq k \leq r$,
    let
    \begin{equation*}
      \calt_n(k) :=
      \Bigl\{T=(T_1,\dots,T_s) \in \calt_n
      : \bigl|{\textstyle\bigcup_{i \in [s]} T_i}\bigr| = k\Bigr\}.
    \end{equation*}
    Note that $\calt_n=\bigcup_{0\leq k\leq r}\calt_n(k)$.
    Observe that $\abs{\calt_n(k)}\leq\binom{\binom{n}{2}}{k}(2^s-1)^k$.
    Indeed, there are $\binom{\binom{n}{2}}{k}$ ways of choosing~$k$
    edges from~$E(K_n)$, and $(2^s-1)^k$ ways of assigning each of
    those edges to the sets~$T_i$ of the $s$-tuples
    $T=(T_1,\dots,T_s)$.  Therefore,
    \begin{equation}
        \label{ineq:sum_of_ps}
        \sum_{T \in \calt_n}
        q^{\big|\bigcup_{i \in [s]} T_i\big|}
        =
        \sum_{k=0}^r
        \abs{\calt_n(k)}
        q^k
        \leq
        \sum_{k=0}^r
        \binom{\binom{n}{2}}{k}
        (2^s)^k
        q^k
        \leq
        1+
        \sum_{k=1}^r
        \(\frac{\ee 2^{s-1} n^2 q}{k}\)^k.
    \end{equation}
    Let $b = 2^{s-1}n^2q$ and $f(x)=(\ee b/x)^x$ for all~$x>0$.  Observe
    that~$f$ is unimodal and achieves its maximum at~$x=b$.  Moreover,
    by~(\ref{eq:A_cond}) we obtain
    \begin{equation*}
      r=\lfloor
      cn^{2-1/m_2(\Harr)}\rfloor\leq2^{s-1}An^{2-1/m_2(\Harr)}
      =2^{s-1}n^2q=b.
    \end{equation*}
    Thus, by~\eqref{ineq:sum_of_ps},
    \begin{equation}
      \sum_{T \in \calt_n}
      q^{\big|\bigcup_{i \in [s]}T_i\big|}
      \leq1+r\left(\frac{\ee b}{r}\right)^r
      =1+r\left(\ee2^{s-1}n^2q\over\bigl\lfloor cn^{2-1/m_2(\Harr)}\bigr\rfloor\right)^r
      \leq2r\left(\ee2^sA\over c\right)^r.
    \end{equation}
    Recalling~\eqref{eq:sb}, we obtain
    \begin{multline}
      \label{eq:3}
      \qquad
      \PP\bigl(G(n,q) \in \calg\bigr)
      \leq2^{r+1}r\ee^{-\alpha n^2q}\left({\ee2^sA\over c}\right)^r
      \leq \ee^{-\alpha n^2q}\left({\ee2^{s+2}A\over c}\right)^r\\
      \leq A^r\ee^{-\alpha n^2q}
      \leq\Bigl(A^{cn^{2-1/m_2(\Harr)}/n^2q}\ee^{-\alpha}\Bigr)^{n^2q}
      =\Bigl(A^{c/A}\ee^{-\alpha}\Bigr)^{n^2q}.\qquad
    \end{multline}
    Owing to~\eqref{eq:A_cond}, we have~$A^{c/A}\ee^{-\alpha}\leq
    \ee^{-\alpha/2}$.  Thus, \eqref{eq:3}~tells us that
    \begin{equation*}
      \PP\bigl(G(n,q) \in \calg\bigr)
      \leq \ee^{-(\alpha/2)n^2q}=e^{-\gamma n^2q},
    \end{equation*}
    as promised in~\eqref{e:exposition}.

    To conclude the proof, we show, using a standard multiple
    exposition argument, that Lemma~\ref{lem:one-statement_exp}
    follows from~\eqref{e:exposition} with $B=2A$ and
    $\beta = \gamma/2$.  Let $p\ge Bn^{-1/m_2(\Harr)} = 2q$.  Note
    that $G\narr \Harr$ is a decreasing property in~$G$.  Let
    $t=\lfloor p/q\rfloor\ge \lfloor B/A \rfloor = 2$ and let
    $G=G_1\cup\dots \cup G_t$, where each~$G_i$ is an independent copy
    of~$G(n,q)$.  It is a simple fact that~$G$ is a~$G(n,p')$ with
    $p'\leq p$.  Also, since $p/q \ge 2$, we have that
    $t \ge p/q-1 \ge p/2q$.  Thus
  \begin{align*}
    \PP(G(n,p)\narr \Harr)
    \stackrel{(\dagger)}{\le}
    \PP(G\narr\Harr)
    \stackrel{(\dagger)}{\le}
            \PP\left(\bigcap_{i\in[t]} \{G_i\narr\Harr\}\right)
     \stackrel{(\sharp)}{=}
            \prod_{i\in[t]} \PP(G_i\narr\Harr)
    & \stackrel{\clap{\eqref{e:exposition}}}{\le}
            \exp(-t\alpha qn^2)                   \\
    & \le   \exp(-\alpha pn^2/2)                  \\
    & =     \exp(-\beta pn^2),
  \end{align*}
  where $(\dagger)$ follow by monotonicity and
  $(\sharp)$ follows from the independence of the $G_i$.
  This concludes the proof of Lemma~\ref{lem:one-statement_exp}.
\end{proof}

The strategy for proving Theorem~\ref{t:small-xmas} is simple: we find
many vertex-disjoint copies of~$\ttt$, and then locate a copy
of~$\Farr$ in the subgraph induced by the roots of those copies.  For
the first step, since $K_4 \arr \ttt$, it suffices to find many
vertex-disjoint copies of~$K_4$.

\begin{lemma}\label{l:number-K4-intro}
  If~$0<\delta \le 1/6$
  and~$p\geq n^{\delta-1/m(K_4)}=n^{\delta-2/3}$, then a.a.s.~$G(n,p)$
  contains at least $cn^{6\delta}$ vertex-disjoint copies of~$K_4$,
  where~$c>0$ is an absolute constant.
\end{lemma}

Lemma~\ref{l:number-K4-intro} is a particular case of Theorem~4
in~\cite{kreuter96:_thres_ramsey} (see also {\cite[Theorem
  3.29]{janson00:_random_graph}}).  For the second step in the proof
of Theorem~\ref{t:small-xmas}, we apply
Lemma~\ref{lem:one-statement_exp}.

\begin{proof}[Proof of Theorem~\ref{t:small-xmas}]
  Let~$\Farr$ be as in the statement of the theorem.  By
  Proposition~\ref{p:m2<2}\,\ref{i:F_a}, we have
  $m_2(\Farr\circ\ttt')=2$.  We will show that there exists
  \begin{equation}
    \label{eq:5}
    p\ll n^{-1/2}=n^{-1/m_2(\Farr\circ\ttt')}
  \end{equation}
  such that a.a.s.\ we have $G(n,p)\arr \Farr\circ\ttt'$.
  Since~$m_2(\Farr)<2$, we can fix positive constants~$\gamma$
  and~$\epsilon$ such that
  \begin{equation}
    \label{eq:gamma_epsilon_def}
    1-\gamma={1\over2}m_2(\Farr)+\epsilon>{1\over2}m_2(\Farr).
  \end{equation}
  Let
  \begin{equation}
    \label{eq:xmas_delta_def}
    \delta={1\over6}\left(1-{6\epsilon\over m_2(\Farr)}\).
  \end{equation}
  By choosing~$\epsilon$ small enough, we may suppose that
  \begin{equation}
    \label{eq:6delta_cond}
    6\delta\geq1-\gamma.
  \end{equation}
  Let~$B$ and~$\beta$ be the constants given by
  Lemma~\ref{lem:one-statement_exp} for~$\Farr$, and let~$c$ be the
  constant in Lemma~\ref{l:number-K4-intro}.  Let
  $B'=c^{-1/m_2(\Farr)}B$.  We may assume that~$B'\geq1$.
  Let~$p=p(n)=B'n^{-(1-\gamma)/m_2(\Farr)}$ and note that,
  by~\eqref{eq:gamma_epsilon_def}, relation~\eqref{eq:5} holds.  We
  show that this choice of~$p$ will do.

  We first claim that, a.a.s., $G=G(n,p)$~is such that
  \begin{itemize}
  \item[(A)] \textsl{any~$U\subseteq V(G)$ with~$|U|\geq cn^{6\delta}$ is
      such that $G[U]\to\Farr$.}
  \end{itemize}
  Indeed, let $t=cn^{1-\gamma}$ and note that
  $p=B'n^{-(1-\gamma)/m_2(\Farr)}=B'c^{1/m_2(\Farr)}t^{-1/m_2(\Farr)}=Bt^{-1/m_2(\Farr)}$.
  By our choice of~$B$ and~$\beta$, for each~$W\subset V(G)$ with
  $|W|=t$, we have
  \begin{equation}
    \label{eq:8}
    \PP\bigl(G[W]\not\arr\Harr\bigr)\leq \ee^{-\beta pt^2}.
  \end{equation}
  Moreover, the number of~$W\subset V(G)$ with~$|W|=t$ is
  ${n\choose t}\leq n^t$.  Since~$m_2(\Farr)>1$, we see that
  \begin{equation}
    \label{eq:1}
    \beta pt^2=\beta tB'n^{-(1-\gamma)/m_2(\Farr)}cn^{1-\gamma}
    =\beta B'ctn^{(1-\gamma)(1-1/m_2(\Farr))}
    \gg t\log n.
  \end{equation}
  Therefore, \eqref{eq:8}~and~\eqref{eq:1} and the union bound tell us
  that, a.a.s., for every~$W\subset V(G)$ with~$|W|=t$, we
  have~$G[W]\to\Farr$.  This implies that $G[U]\to\Farr$ for every~$U$
  with~$|U|\geq cn^{6\delta}$, because $t=cn^{1-\gamma}$ and we
  have~\eqref{eq:6delta_cond}, concluding the proof of claim~(A).

  Now note that, by~\eqref{eq:gamma_epsilon_def}, we have
  $p=B'n^{-(1-\gamma)/m_2(\Farr)}\geq n^{-(1-\gamma)/m_2(\Farr)}
  =n^{-1/2-\epsilon/m_2(\Farr)}$.  By~\eqref{eq:xmas_delta_def}, we
  have $-1/2-\epsilon/m_2(\Farr)=\delta-2/3=\delta-1/m(K_4)$.  Thus,
  Lemma~\ref{l:number-K4-intro} tells us that, a.a.s.,
  \begin{itemize}
  \item[(B)] \textsl{$G$~contains at least~$cn^{6\delta}$
      vertex-disjoint copies of~$K_4$.}
  \end{itemize}

  Suppose~$G$ satisfies~(A) and~(B).  We show that
  $G\to\Farr\circ\ttt'$.  Since~(A) and~(B) hold a.a.s., this will
  conclude the proof.  Let~$\Garr$ be an orientation of~$G$.  Noticing
  that $K_4\arr\ttt$, property~(B) implies that~$\Garr$
  contains~$cn^{6\delta}$ vertex-disjoint copies of~$\ttt'$.
  Let~$U\subset V(G)$ be the set of~$cn^{6\delta}$ roots of those
  copies of~$\ttt'$.  Since~(A) holds, there is
  a copy of~$\Farr$ in~$\Garr[U]$.
  Therefore,~$\Garr$ contains a copy
  of~$\Farr\circ\ttt'$.  This shows that $G\to\Farr\circ\ttt'$, as
  required.
\end{proof}

\section{Rooted products with larger trees}
\label{s:large-xmas-tree}

We now prove Theorem~\ref{t:large-xmas-tree}.  We shall use the
following corollary of~\cite[Theorem 3]{naia22}.

\begin{theorem}[{\cite{naia22}}]\label{t:Burr-with-logn-factor}
  If $G$ is a graph and~$\Tarr$ is an oriented tree of
  order~$\lfloor\chi(G)/\lceil\log_2 v(G)\rceil\rfloor$,
  then $G\arr \Tarr$.
\end{theorem}

\begin{proof}[Proof of Theorem~\ref{t:large-xmas-tree}]
  Let $b=c(\log2)/4$, where~$c$ is the constant in
  Lemma~\ref{l:number-K4-intro}.
  Recall that $p=n^{\delta-2/3}$ and that~$2/21<\delta\le 1/6$.
  Let~$\Farr$ be an oriented tree with
  $v(\Farr)\le bn^{7\delta-2/3}/(\log n)^2$.  As in the proof of
  Theorem~\ref{t:small-xmas}, it suffices to prove that,
  for $p=n^{\delta-2/3}$ as in the statement
  of~Theorem~\ref{t:large-xmas-tree},
  a.a.s.,
  $G=G(n,p)$~has the following two properties:
  \begin{itemize}
  \item[(A)]$G[U]\to\Farr$ for any $U\subseteq V(G)$ with~$|U|\geq
    cn^{6\delta}$.
  \item[(B)]$G$ contains at least~$cn^{6\delta}$ vertex-disjoint
    copies of~$K_4$.
  \end{itemize}
  The fact that~(B) holds a.a.s.\ is Lemma~\ref{l:number-K4-intro}.
  We derive that~(A) holds a.a.s.\ from
  Theorem~\ref{t:Burr-with-logn-factor}.  An easy application of the
  first moment method shows that a.a.s.\ $\alpha(G) < 3p^{-1}\log n$,
  and hence we suppose that this inequality holds for~$G$.
  Fix~$U\subseteq V(G)$ with~$|U|\geq cn^{6\delta}$.  Then
  \begin{equation}\label{e:chi(G[S])}
    \chi(G[U])
    \ge\frac{|U|}{\alpha(G(n,p))}
    >\frac{cn^{6\delta}p}{3\log n}
    =\frac{cn^{7\delta-2/3}}{3\log n}.
  \end{equation}
  Thus, for any large enough~$n$,
  \begin{equation*}
    \left\lfloor\frac{\chi(G[U])}{\lceil\log_2|U|\rceil}\right\rfloor
    \stackrel{\eqref{e:chi(G[S])}}{>}
    \frac{cn^{7\delta-2/3}}{3(\log n)(\log_2n+1)}-1
    \geq\left(\frac{c\log2}4\right){n^{7\delta-2/3}\over(\log n)^2}
    ={bn^{7\delta-2/3}\over(\log n)^2}.
  \end{equation*}
  Therefore~$v(\Farr)\leq{bn^{7\delta-2/3}/(\log n)^2}
  \leq\chi(G[U])/\lceil\log_2(G[U])\rceil$ and
  Theorem~\ref{t:Burr-with-logn-factor} tells us that, indeed,
  $G[U]\to\Farr$.
\end{proof}

\section{Concluding remarks}
\label{sec:concluding-remarks}

We believe that finding the threshold~$p_\Harr$ for the
property~$G(n,p)\to\Harr$ for a given acyclically oriented graph~$\Harr$
is a natural problem.  Methods developed in this general line of
research show that~$p_\Harr\leq n^{-1/m_2(\Harr)}$, but we have shown
here that we may have~$p_\Harr\ll n^{-1/m_2(\Harr)}$ for
certain~$\Harr$.  It would be interesting to characterise
those~$\Harr$ for which~$p_\Harr=n^{-1/m_2(\Harr)}$.

A related problem is establishing sharp thresholds for
the appearance of a given fixed-size oriented tree~$\Tarr$.
Indeed, it can be shown that~$p_{\Tarr}=\Theta(1/n)$
(this follows from a theorem of Burr~\cite{Burr80}:
 if $\Tarr$ is
an oriented tree of order $t$ and $G$ is a graph with
chromatic number $(t-1)^2$, then $G\arr\Tarr$). And yet,
determining the precise constant for an arbitrary tree
might be challenging, possibly even for paths.
We remark that the orientation of $\Tarr$ plays a role
in the threshold, as distinct orientations
of a $t$-vertex path can be shown to have distinct
thresholds
(this shall be discussed in future work by a subset of the authors).

\section*{Acknowledgements}

We thank the anonymous referees for their careful reading and helpful
comments.

\begin{adjustwidth}{0em}{-5.05em}
  \bibliographystyle{abbrv}
{\footnotesize
  \bibliography{extracted-mrefed.bib}
}
\end{adjustwidth}

\end{document}